\documentclass[preprint,10pt]{elsarticle}

\makeatletter
\def\ps@pprintTitle{%
   \let\@oddhead\@empty
   \let\@evenhead\@empty
   \def\@oddfoot{\reset@font\hfil\thepage\hfil}
   \let\@evenfoot\@oddfoot
}
\makeatother

\usepackage[margin=1in]{geometry}
\usepackage{amsmath,amsfonts,mathtools,amssymb}
\usepackage[T1]{fontenc}\usepackage{lmodern}
\usepackage{booktabs}
\usepackage{subcaption}
\usepackage{xcolor}
\usepackage{hyperref}
\hypersetup{hidelinks}
\usepackage{multirow}
\usepackage{placeins}
\usepackage{etoolbox}
\usepackage{doi}

\usepackage{tikz}
\usetikzlibrary{calc}
\usepackage{tkz-euclide}
\usetkzobj{all} 
\usepackage{calc}
\usepackage{tikz,pgfplots}
\usetikzlibrary{arrows.meta,3d,matrix,positioning,decorations.markings}

\DeclareMathOperator*{\argmin}{arg\,min}

\newcommand{\maxpoly}{d}

\newcommand{\ngroup}{{n_g}}
\newcommand{\groupindex}{k}

\newcommand{\nfunc}{{n_\func}}
\newcommand{\nsing}{{n_\singfunc}}
\newcommand{\nmono}{{n_\monofunc}}

\newcommand{\func}{f}
\newcommand{\singfunc}{s}
\newcommand{\monofunc}{m}

\newcommand{\funcvec}{\mathbf{\func}}
\newcommand{\singfuncvec}{\mathbf{\singfunc}}
\newcommand{\monofuncvec}{\mathbf{\monofunc}}

\newcommand{\npoints}{n}

\newcommand{\pz}{\phantom{0}}

\definecolor{orange}{rgb}{1,0.5,0}
\definecolor{green}{rgb}{0,0.5,0}
\newcommand{\reviewerOne}[1]{{#1}}
\newcommand{\reviewerTwo}[1]{{#1}}
\begin{document}

\begin{frontmatter}



\title{Symmetric Triangle Quadrature Rules for Arbitrary Functions}


\author{Brian A.\ Freno}
\ead{bafreno@sandia.gov}
\author{William A.\ Johnson}
\ead{wajohns@sandia.gov}
\author{Brian F.\ Zinser}
\ead{bzinser@sandia.gov}
\author{Salvatore Campione}
\ead{sncampi@sandia.gov}
\address{Sandia National Laboratories, PO Box 5800, Albuquerque, NM 87185}

\begin{abstract}
Despite extensive research on symmetric polynomial quadrature rules for triangles, as well as approaches to their calculation, few studies have focused on non-polynomial functions, particularly on their integration using symmetric triangle rules.  In this paper, we present two approaches to computing symmetric triangle rules for singular integrands by developing rules that can integrate arbitrary functions.  The first approach is well suited for a moderate amount of points and retains much of the efficiency of polynomial quadrature rules.  The second approach better addresses large amounts of points, though it is less efficient than the first approach.  We demonstrate the effectiveness of both approaches on singular integrands, which can often yield relative errors two orders of magnitude less than those from polynomial quadrature rules.
\end{abstract}


\end{frontmatter}


\section{Introduction}

Due to their efficiency, Gaussian quadrature rules are useful for numerical integration.  For integrands that can be accurately approximated by polynomials, rules are typically employed that exactly integrate polynomials of increasing degree.  

Because of their common use in two-dimensional discretizations, the development of quadrature rules for triangles is a popular research area.  Several authors have developed methods for computing symmetric quadrature rules for polynomials~\cite{lyness_1975,dunavant_1985,wandzura_2003,papanicolopulos_2015}.  Symmetric rules are desirable because their mapping to the integration domain is straightforward and points are not heavily concentrated near some vertices.  Asymmetric rules, on the other hand, require the determination of the vertex mapping, and point concentration at the vertices is inconsistent.

In Reference~\cite{lyness_1975}, the authors present quadrature rules for many polynomial degrees, up to degree 12.  In Reference~\cite{dunavant_1985}, the author provides quadrature rules for all degrees, up to degree 20.  Reference~\cite{wandzura_2003} uses numerical optimization to compute even higher degrees of polynomials.  In two dimensions, the optimal number of integration points is much less straightforward than for one dimension, and, for a given degree and number of points, there can be multiple solutions for the points and weights.  Reference~\cite{papanicolopulos_2015} presents an approach for determining all of the solutions for polynomials up to moderately high degrees, and the author presents these solutions as ancillary material in Reference~\cite{papanicolopulos_2015_anc}.  Reference~\cite{mousavi_2010} presents an approach for computing quadrature rules for polynomials over arbitrary polygons.

However, for functions with singularities on the integration boundary, rules for polynomials do not converge monotonically or as rapidly as the number of integration points is increased.  Provided they are integrable, these singular integrands may include unbounded derivatives at the boundary, where the integrand may be defined or undefined.  In Reference~\cite{ma_1996}, the authors present an approach for computing the quadrature rules associated with arbitrary one-dimensional functions and demonstrate the effectiveness of their approach for several sequences of functions, with and without singularities.  To extend that approach to two dimensions, the authors of Reference~\cite{vipiana_2013} take an outer product of the one-dimensional rules and asymmetrically map the result to a triangle.

Regardless of dimension and function sequence, the equations for computing quadrature rules are stiff and highly dependent on the initial guess~\cite{dunavant_1985,ma_1996,wandzura_2003,papanicolopulos_2015}.  For one dimension, the authors of Reference~\cite{ma_1996} use a continuation method, beginning with the polynomial quadrature rules and gradually transforming the sequence of polynomials to the desired sequence of functions, using the intermediate solutions as subsequent initial guesses.  In higher dimensions, the problem is complicated by the potentially unknown number of optimal points~\cite{xiao_2010} and the potential existence of multiple solutions~\cite{papanicolopulos_2015,papanicolopulos_2015_anc}. 

In this paper, we develop symmetric quadrature rules for triangles that integrate arbitrary functions, motivated by the need to integrate integrands with boundary singularities.  This paper is organized as follows.  In Section~\ref{sec:preliminaries}, we discuss the details of symmetric quadrature rules for triangles.  In Section~\ref{sec:singularities}, we discuss singularities and the construction of one- and two-dimensional function sequences to be integrated exactly.  In Section~\ref{sec:moderate}, we describe our first approach to computing symmetric quadrature rules, which is better suited for moderate amounts of functions and points (leading to about six or seven digits in accuracy), and, in Section~\ref{sec:large}, we describe our second approach to computing symmetric quadrature rules, which is better suited for large amounts of functions and points  (leading to machine accuracy).  In Section~\ref{sec:example}, we demonstrate the two proposed approaches for a sample triangle and compare with polynomial rules.  Finally, in Section~\ref{sec:conclusions}, we provide an outlook for future work.

\section{Quadrature Preliminaries}
\label{sec:preliminaries}

In this section, we describe the concepts we use to construct our approaches to computing symmetric quadrature rules for triangles that integrate singularities.

\subsection{Quadrature Rules} 
An $\npoints$-point quadrature rule exactly integrates a sequence of $\nfunc$ functions $\mathbf{f}(\mathbf{x})=\{f_1(\mathbf{x}),\hdots,f_\nfunc(\mathbf{x})\}$, such that
\begin{align}
\int_A \mathbf{f}(\mathbf{x})dA = \sum_{i=1}^\npoints w_i \mathbf{f}(\mathbf{x}_i).
\label{eq:quadrature}
\end{align}
In Equation~\eqref{eq:quadrature}, the quadrature integration exactly computes the integrals by taking a linear combination of the function values at $\mathbf{x}_i$, which are weighted by weights $w_i$, for $i=1,\hdots,\npoints$.
In one dimension, $\nfunc=2n$ and, for polynomials, $\mathbf{f}(x)=\{1,\hdots,x^{2n-1}\}$.  For two dimensions, one could speculate $\nfunc=3n$~\cite{xiao_2010}, but the ability to achieve such efficiency is unproven, and, if the rules are required to be symmetric, the efficiency can be significantly lower~\cite{xiao_2010}.

\subsection{Symmetric Rules for Triangles} 
\label{sec:symmetry}

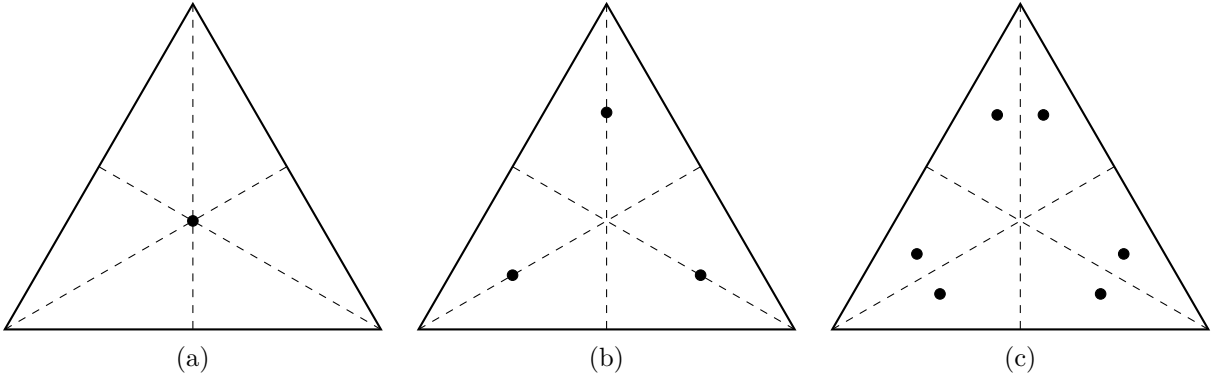
\begin{figure}
\centering
\begin{tikzpicture}

\def\ts{5}; 

\coordinate (A) at (-0.5,{-sqrt(3)/6});
\coordinate (B) at ( 0.5,{-sqrt(3)/6});
\coordinate (C) at ( 0  , {sqrt(1/3)});
\coordinate (D) at ($0.5*(A)+0.5*(B)$);
\coordinate (E) at ($0.5*(B)+0.5*(C)$);
\coordinate (F) at ($0.5*(C)+0.5*(A)$);
\coordinate (O) at (0,0);

\def\tx{0}
\def\ty{{\ts*sqrt(3)/6}}

\coordinate (T1) at ($\ts*(A)+(\tx,\ty)$);
\coordinate (T2) at ($\ts*(B)+(\tx,\ty)$);
\coordinate (T3) at ($\ts*(C)+(\tx,\ty)$);
\coordinate (T4) at ($\ts*(D)+(\tx,\ty)$);
\draw[thick] (T1) -- (T2) -- (T3) -- cycle;

\draw[dashed] ($\ts*(C)+(\tx,\ty)$) -- ($\ts*(D)+(\tx,\ty)$);
\draw[dashed] ($\ts*(A)+(\tx,\ty)$) -- ($\ts*(E)+(\tx,\ty)$);
\draw[dashed] ($\ts*(B)+(\tx,\ty)$) -- ($\ts*(F)+(\tx,\ty)$);

\draw[fill=black] ($\ts*(0,0)+(\tx,\ty)$) circle (.07);

\node[below = .1 of T4] {(a)};

\def\tx{5.5}
\def\ty{{\ts*sqrt(3)/6}}

\coordinate (T1) at ($\ts*(A)+(\tx,\ty)$);
\coordinate (T2) at ($\ts*(B)+(\tx,\ty)$);
\coordinate (T3) at ($\ts*(C)+(\tx,\ty)$);
\coordinate (T4) at ($\ts*(D)+(\tx,\ty)$);
\draw[thick] (T1) -- (T2) -- (T3) -- cycle;

\draw[dashed] ($\ts*(C)+(\tx,\ty)$) -- ($\ts*(D)+(\tx,\ty)$);
\draw[dashed] ($\ts*(A)+(\tx,\ty)$) -- ($\ts*(E)+(\tx,\ty)$);
\draw[dashed] ($\ts*(B)+(\tx,\ty)$) -- ($\ts*(F)+(\tx,\ty)$);

\draw[fill=black] ($\ts*(-0.25,-0.144337567297406)+(\tx,\ty)$) circle (.07);
\draw[fill=black] ($\ts*( 0.25,-0.144337567297406)+(\tx,\ty)$) circle (.07);
\draw[fill=black] ($\ts*( 0.00, 0.288675134594813)+(\tx,\ty)$) circle (.07);

\node[below = .1 of T4] {(b)};


\def\tx{11}
\def\ty{{\ts*sqrt(3)/6}}

\coordinate (T1) at ($\ts*(A)+(\tx,\ty)$);
\coordinate (T2) at ($\ts*(B)+(\tx,\ty)$);
\coordinate (T3) at ($\ts*(C)+(\tx,\ty)$);
\coordinate (T4) at ($\ts*(D)+(\tx,\ty)$);
\draw[thick] (T1) -- (T2) -- (T3) -- cycle;
\draw[dashed] ($\ts*(C)+(\tx,\ty)$) -- ($\ts*(D)+(\tx,\ty)$);
\draw[dashed] ($\ts*(A)+(\tx,\ty)$) -- ($\ts*(E)+(\tx,\ty)$);
\draw[dashed] ($\ts*(B)+(\tx,\ty)$) -- ($\ts*(F)+(\tx,\ty)$);

\draw[fill=black] ($\ts*( 0.061447179740077, 0.282059528176809)+(\tx,\ty)$) circle (.07);
\draw[fill=black] ($\ts*(-0.274994306650607,-0.087814945442589)+(\tx,\ty)$) circle (.07);
\draw[fill=black] ($\ts*( 0.213547126910531,-0.194244582734219)+(\tx,\ty)$) circle (.07);
\draw[fill=black] ($\ts*( 0.274994306650607,-0.087814945442589)+(\tx,\ty)$) circle (.07);
\draw[fill=black] ($\ts*(-0.061447179740077, 0.282059528176809)+(\tx,\ty)$) circle (.07);
\draw[fill=black] ($\ts*(-0.213547126910531,-0.194244582734219)+(\tx,\ty)$) circle (.07);

\node[below = .1 of T4] {(c)};

\end{tikzpicture}
\caption{Pictorial representation of (a) type-0, (b) type-1, and (c) type-2 orbits.}
\label{fig:three_orbits}
\end{figure}

Symmetric rules for triangles are invariant to rotation and reflection about the medians for equilateral triangles, which can be isoparametrically transformed to arbitrary triangles.  As several references~\cite{lyness_1975,dunavant_1985,wandzura_2003,papanicolopulos_2015} have described, the points are comprised of a combination of orbits.  There are three types of orbits, which are shown in Figure~\ref{fig:three_orbits}.

The type-0 orbit consists of a point at the centroid, which is located at $\left(\frac{1}{3},\frac{1}{3},\frac{1}{3}\right)$ in barycentric coordinates.  The type-1 orbit consists of three points, each on a median, such that the coordinates are the three unique permutations of $\left(\alpha,\frac{1-\alpha}{2},\frac{1-\alpha}{2}\right)$.  The type-2 orbit consists of six points, not on the medians, such that the coordinates are the six unique permutations of $(\alpha,\beta,1-\alpha-\beta)$.

In terms of orbits, the number of points $\npoints$ is
\begin{align*}
\npoints = n_0 + 3n_1 + 6n_2,
\end{align*}
where $n_j$ is the number of type-$j$ orbits.  For type-0 orbits, $n_0=0$ or $n_0=1$.  Type-1 and type-2 orbits can have arbitrary $n_1$ and $n_2$.  We denote the orbit counts for a given $\npoints$ as the triplet $(n_0,n_1,n_2)$.  For each orbit, the weights $w_i$ for the points are the same.  A type-0 orbit has an unknown weight, a type-1 orbit has an unknown weight and coordinate, and a type-2 orbit has an unknown weight and two unknown coordinates.  Therefore, the total number of unknowns is $n_0+2n_1+3n_2$.

\subsection{Polynomial Integration} 
\label{sec:polynomials}

Polynomial rules capable of integrating polynomials up to degree $\maxpoly$ can exactly integrate linear combinations of the monomials $x^p y^q$, where $0\le p\le \maxpoly$, $0\le q\le \maxpoly$, $0\le p+q\le \maxpoly$, totaling $\nfunc=(\maxpoly+1)(\maxpoly+2)/2$ monomials.  This requirement can yield more equations than unknowns.  For example, a $(0,1,0)$ rule can integrate polynomials up to $\maxpoly=2$, resulting in $\nfunc=6$ monomials: $\mathbf{f}(x,y)=\{1,\,x,\,y,\,x^2,\,y^2,\,xy\}$.  However, the number of unknowns is two: $\alpha$ and $w$.  This mismatch is reconcilable because the $\nfunc=6$ equations of Equation~\eqref{eq:quadrature} are not linearly independent.

To reduce the number of equations to an amount that is linearly independent, one can use de Moivre's theorem to express the polynomials in polar coordinates~\cite{lyness_1975,dunavant_1985} or otherwise construct a sequence of invariant polynomials~\cite{wandzura_2003,papanicolopulos_2015}.  

On the other hand, instead of formulating the problem as the solution to a system of equations, one can formulate the problem as an unconstrained optimization problem:
\begin{align}
\argmin_{\boldsymbol{\alpha},\boldsymbol{\beta},\mathbf{w}} F(\boldsymbol{\alpha},\boldsymbol{\beta},\mathbf{w}),
\label{eq:opt_bary}
\end{align}
where
\begin{align}
\reviewerOne{%
F(\boldsymbol{\alpha},\boldsymbol{\beta},\mathbf{w}) = \sum_{j=1}^\nfunc \left(\frac{\tilde{I}_{f_j}-I_{f_j}}{I_{f_j}}\right)^2,
\label{eq:objective_function}
}
\end{align}
with the expectation that the objective function $F(\boldsymbol{\alpha},\boldsymbol{\beta},\mathbf{w})$ is reduced to zero.  
\reviewerOne{%
In Equation~\eqref{eq:objective_function}, 
$\tilde{I}_{f_j}=\sum_{i=1}^\npoints w_i' f_j(\alpha_i,\beta_i)$
and
$I_{f_j}=\int_0^1\int_0^{1-\beta} f_j(\alpha,\beta)d\alpha d\beta$.}
\reviewerTwo{The scaling provided by the denominator in Equation~\eqref{eq:objective_function} can improve optimization convergence but should be omitted if any of the integrals are zero.}  Additionally, the integrals in Equation~\eqref{eq:objective_function} should be computed analytically when possible.
In Equation~\eqref{eq:objective_function}, because of symmetry, the ability to integrate $f(\alpha,\beta)$ indicates the ability to integrate $f(\beta,\alpha)$, enabling the number of functions to be reduced to $\nfunc=\left\lfloor\frac{(d+2)^2}{4}\right\rfloor$, where $\lfloor\cdot\rfloor$ denotes the floor function.  For example, when $d=3$, the $\nfunc=6$ functions are $\mathbf{f}(\alpha,\beta)=\{1,\,\alpha,\,\alpha^2,\,\alpha\beta,\,\alpha^3,\,\alpha^2\beta\}$.

\section{Singularities and Function-Sequence Construction} 
\label{sec:singularities}
Integrands with boundary singularities can have singularities located on edges and/or vertices.  At these locations, derivatives of the integrands are unbounded, and the integrand can be defined or undefined, provided the integrand is integrable.  For example, in electromagnetic simulations, logarithmic singularities can appear at the integration boundaries~\cite{vipiana_2013,freno_ieee}.  We present two approaches for constructing the sequence of functions to be exactly integrated by the quadrature rules.

\subsection{One-Dimensional Functions}
\label{sec:singularity_1d}

The integrands can be analyzed in terms of series expansions about the edges or vertices.  These expansions can take the form of alternating monomials and singularities; for example, $\funcvec(x)=\{1,\allowbreak\,\ln x,\allowbreak\,x,\allowbreak\, x\ln x,\hdots\}$~\cite{ma_1996} or $\funcvec(x)=\{1,\allowbreak\,x,\allowbreak\,x\ln x,\allowbreak\,x^2,\allowbreak\,x^3,\allowbreak\,x^3\ln x,\allowbreak\,x^4,\allowbreak\,x^5,\allowbreak\,x^5\ln x,\hdots\}$~\cite{freno_ieee}.  Note that, given the electromagnetic applications in mind, we concentrate on logarithmic functions here, but other singular functions can be considered.  More generally, we denote these expansions as $\funcvec(x) = \{\monofuncvec(x),\,\singfuncvec(x)\}$, where
$\monofuncvec(x)=\{1,\,x,\hdots,x^\maxpoly\}$, and
$\singfuncvec(x)=\{s_1(x),\allowbreak\,s_2(x),\hdots,s_\nsing(x)\}$ is the sequence of singular functions.  \reviewerOne{For a given maximum polynomial degree and number of singularities, the} total number of functions is $\nfunc=\nmono+\nsing$ and the number of monomials is $\nmono=\maxpoly+1$.

For triangles, in barycentric coordinates, edge singularities ($\alpha=0$) can be modeled by $\singfuncvec(\alpha)$ and vertex singularities ($\alpha=1$) can be modeled by $\singfuncvec(1-\alpha)$, assuming the singularities in $\singfuncvec(x)$ occur at $x=0$.  For both vertex and edge singularities, one can increase the number of entries in $\singfuncvec(x)$; however, doing so may further reduce the maximum polynomial degree the rules are capable of integrating for a given number of points.

\subsection{Two-Dimensional Functions}
\label{sec:singularity_2d}

For this approach, $\monofuncvec(x,y)$ is the sequence of $\nmono=\left\lfloor\frac{(d+2)^2}{4}\right\rfloor$ monomials described in Section~\ref{sec:polynomials}.
\reviewerOne{%
For example, for $d=3$, $\nmono=6$ and $\monofuncvec(x,y)=\{1,\allowbreak\,x,\allowbreak\,x^2,\allowbreak\,xy,\allowbreak\,x^3,\allowbreak\,x^2y\}$.}
If a two-dimensional characterization of the singularities is known, it can be included in the singularity sequence $\singfuncvec(x,y)$; otherwise, the singularities can be obtained from series expansions, as is done in Section~\ref{sec:singularity_1d} for $\singfuncvec(x)$.  The total number of functions is $\nfunc=\nmono+\nsing$.  

We introduce the concept of function groups as an accounting mechanism\reviewerTwo{, denoting the function group index by $\groupindex$ and the highest function group index by $\ngroup$, such that $\groupindex=0,\hdots,\ngroup$}.  
A function group \reviewerOne{$\mathbf{f}^\groupindex(x,y)$} contains either one singularity or the monomials with powers that sum to a particular degree.  \reviewerOne{In the latter case, the function group associated with a polynomial of degree $\maxpoly'$ contains $1+\lfloor \maxpoly'/2\rfloor$ monomials.  
In order to integrate a given function group, we require that the quadrature rules be able to integrate the preceding groups.  Therefore, the function sequence is $\mathbf{f}(x,y)=\mathbf{f}_{\ngroup}(x,y)=\{\mathbf{f}^0(x,y),\hdots,\mathbf{f}^\ngroup(x,y)\}$.}

\section{Approach 1: Optimization for a Moderate Number of Functions}
\label{sec:moderate}

The goal of this work is to achieve the ability to integrate polynomials \reviewerOne{with an efficiency similar to that of the} symmetric polynomial rules, while being able to integrate singularities.  Therefore, this approach uses the polynomial rules as a baseline.
\reviewerOne{As an example, for a 13-point $(1,2,1)$ rule, Figure~\ref{fig:opt_fig} depicts the eight unknowns: $\alpha_2$, $\alpha_3$, $\alpha_4$, $\beta_4$, $w_1$, $w_2$, $w_3$, and $w_4$.}

\begin{figure}[h]
\vspace*{-0.5em}
\centering
\begin{tikzpicture}

\definecolor{lightgreen}{RGB}{173.4,217,170.6}
\definecolor{darkgreen} {RGB}{ 51  ,160, 44}
\definecolor{darkblue}  {RGB}{ 31  ,120,180}
\definecolor{lightblue} {RGB}{165.4,201,225}
\definecolor{darkred}   {RGB}{227  ,26,28}
\definecolor{lightred}  {RGB}{243.8,163.4,164.2}
\definecolor{darkorange}{RGB}{255  ,127  ,0}

\def\ts{7.5}; 

\coordinate (A) at (-0.5,{-sqrt(3)/6});
\coordinate (B) at ( 0.5,{-sqrt(3)/6});
\coordinate (C) at ( 0  , {sqrt(1/3)});
\coordinate (D) at ($0.5*(A)+0.5*(B)$);
\coordinate (E) at ($0.5*(B)+0.5*(C)$);
\coordinate (F) at ($0.5*(C)+0.5*(A)$);
\coordinate (O) at (0,0);

\def\tx{0}
\def\ty{{\ts*sqrt(3)/6}}

\coordinate (T1) at ($\ts*(A)+(\tx,\ty)$);
\coordinate (T2) at ($\ts*(B)+(\tx,\ty)$);
\coordinate (T3) at ($\ts*(C)+(\tx,\ty)$);
\coordinate (T4) at ($\ts*(D)+(\tx,\ty)$);

\node[above = 0 of T3] {$(1,0,0)$};
\node[anchor=north] at (T2) {$(0,0,1)$};
\node[anchor=north] at (T1) {$(0,1,0)$};

\draw[dashed] ($\ts*(C)+(\tx,\ty)$) -- ($\ts*(D)+(\tx,\ty)$);
\draw[dashed] ($\ts*(A)+(\tx,\ty)$) -- ($\ts*(E)+(\tx,\ty)$);
\draw[dashed] ($\ts*(B)+(\tx,\ty)$) -- ($\ts*(F)+(\tx,\ty)$);

\draw[draw=darkgreen,fill=darkgreen] ($\ts*(  0.00000000000000,  0.00000000000000)+(\tx,\ty)$) circle (.07);
\draw[draw=darkorange,fill=darkorange] ($\ts*( -0.10948105088144, -0.06320891419756)+(\tx,\ty)$) circle (.07);
\draw[draw=darkorange,fill=darkorange] ($\ts*(  0.10948105088144, -0.06320891419756)+(\tx,\ty)$) circle (.07);
\draw[draw=darkorange,fill=darkorange] ($\ts*(  0.00000000000000,  0.12641782839512)+(\tx,\ty)$) circle (.07);
\draw[draw=darkred,fill=darkred] ($\ts*( -0.40230484564668, -0.23227081093040)+(\tx,\ty)$) circle (.07);
\draw[draw=darkred,fill=darkred] ($\ts*(  0.40230484564668, -0.23227081093040)+(\tx,\ty)$) circle (.07);
\draw[draw=darkred,fill=darkred] ($\ts*(  0.00000000000000,  0.46454162186080)+(\tx,\ty)$) circle (.07);
\draw[draw=darkblue,fill=darkblue] ($\ts*(  0.13208759028978,  0.26423375160519)+(\tx,\ty)$) circle (.07); 
\draw[draw=darkblue,fill=darkblue] ($\ts*(  0.29487693657225, -0.01772566708697)+(\tx,\ty)$) circle (.07); 
\draw[draw=darkblue,fill=darkblue] ($\ts*( -0.13208759028978,  0.26423375160519)+(\tx,\ty)$) circle (.07); 
\draw[draw=darkblue,fill=darkblue] ($\ts*(  0.16278934628247, -0.24650808451821)+(\tx,\ty)$) circle (.07); 
\draw[draw=darkblue,fill=darkblue] ($\ts*( -0.29487693657225, -0.01772566708697)+(\tx,\ty)$) circle (.07); 
\draw[draw=darkblue,fill=darkblue] ($\ts*( -0.16278934628247, -0.24650808451821)+(\tx,\ty)$) circle (.07); 


\node[anchor=north] at ($\ts*(  0.00000000000000,  -0.03000000000000)+(\tx,\ty)$) {$\color{darkgreen}w_1$};

\draw[thick,draw=darkred,fill=darkred,->,>=stealth] ($\ts*(  .05,  {-sqrt(3)/6})+(\tx,\ty)$) -- ($\ts*(  0.05,  0.46454162186080)+(\tx,\ty)$) node [darkred,midway, below,rotate=90] {$\alpha_3$};
\draw[thick,draw=darkred,fill=darkred] ($\ts*(  .04,  0.46454162186080)+(\tx,\ty)$) -- ($\ts*(  0.06,  0.46454162186080)+(\tx,\ty)$);

\node[anchor=north east] at ($\ts*(  0.00000000000000,  0.46454162186080)+(\tx,\ty)$) {$\color{darkred}w_3$};

\draw[thick,draw=darkorange,fill=darkorange,->,>=stealth] ($\ts*(  -.05,  {-sqrt(3)/6})+(\tx,\ty)$) -- ($\ts*(  -0.05,  0.12641782839512)+(\tx,\ty)$) node [darkorange,midway, above,rotate=90] {$\alpha_2$};
\draw[thick,draw=darkorange,fill=darkorange] ($\ts*(  -.04,  0.12641782839512)+(\tx,\ty)$) -- ($\ts*(  -0.06,  0.12641782839512)+(\tx,\ty)$);

\node[anchor=south east] at ($\ts*(  0.00000000000000,  0.12641782839512)+(\tx,\ty)$) {$\color{darkorange}w_2$};

\draw[thick,draw=darkblue,fill=darkblue,->,>=stealth] ($\ts*( 0.16278934628247+.05,  {-sqrt(3)/6})+(\tx,\ty)$) -- ($\ts*(  0.16278934628247+.05,  -0.24650808451821)+(\tx,\ty)$) node [darkblue,near end, below,rotate=90] {$\alpha_4$};
\draw[thick,draw=darkblue,fill=darkblue] ($\ts*(0.16278934628247+.04,  -0.24650808451821)+(\tx,\ty)$) -- ($\ts*(0.16278934628247+.06,  -0.24650808451821)+(\tx,\ty)$);

\draw[thick,draw=darkblue,fill=darkblue,->,>=stealth] ($\ts*( 0.16278934628247, -0.24650808451821)+\ts*.05*(-.5,{sqrt(3)/2})+\ts*0.270949467507843*({sqrt(3)/2},.5)+(\tx,\ty)$) -- ($\ts*( 0.16278934628247, -0.24650808451821)+\ts*.05*(-.5,{sqrt(3)/2})+(\tx,\ty)$) node [darkblue,midway, above,rotate=30] {$\beta_4$}; 
\draw[thick,draw=darkblue,fill=darkblue] ($\ts*( 0.16278934628247, -0.24650808451821)+\ts*.06*(-.5,{sqrt(3)/2})+(\tx,\ty)$) -- ($\ts*( 0.16278934628247, -0.24650808451821)+\ts*.04*(-.5,{sqrt(3)/2})+(\tx,\ty)$);

\node[anchor=east] at ($\ts*(  0.16278934628247, -0.24650808451821)+(\tx,\ty)$) {$\color{darkblue}w_4$};


\draw[thick] (T1) -- (T2) -- (T3) -- cycle;
\end{tikzpicture}
\vspace*{-0.5em}
\caption{Pictorial representation of Approach 1: The eight unknowns computed in Equation~\eqref{eq:opt_bary} for a $(1,2,1)$ rule.}
\label{fig:opt_fig}
\end{figure}

\reviewerOne{For polynomial rules, each function group $\mathbf{f}^k(x,y)$ consists of the monomials with powers that sum to $k$ for $k=1,\hdots,\maxpoly$.  Therefore, for polynomial rules, $\ngroup=\maxpoly$; Table~\ref{tab:polynomial_degrees} lists these values for each number of integration points $\npoints$}, as well as the orbit counts~\cite{dunavant_1985,papanicolopulos_2015}.
%
Because these orbit counts have been shown to be the most efficient~\cite{dunavant_1985,papanicolopulos_2015}, we will use these choices for each $\npoints$ listed in Table~\ref{tab:polynomial_degrees} \reviewerOne{for arbitrary function sequences}.

\begin{table}[htbp!]
\centering
\vspace*{-.5em}
\begin{tabular}{c c c c c | c c c c c}
\toprule
$\npoints$ & $n_0$ & $n_1$ & $n_2$ & $\ngroup=\maxpoly$ & $\npoints$ & $n_0$ & $n_1$ & $n_2$ & $\ngroup=\maxpoly$ \\
\midrule
\pz1 & 1 & 0 & 0 & \pz1              & 27 & 0 & \pz5 & 2 & 11  \\
\pz3 & 0 & 1 & 0 & \pz2              & 33 & 0 & \pz5 & 3 & 12  \\
\pz4 & 1 & 1 & 0 & \pz3              & 37 & 1 & \pz6 & 3 & 13  \\
\pz6 & 0 & 2 & 0 & \pz4              & 42 & 0 & \pz6 & 4 & 14  \\
\pz7 & 1 & 2 & 0 & \pz5              & 48 & 0 & \pz6 & 5 & 15  \\
  12 & 0 & 2 & 1 & \pz6              & 52 & 1 & \pz7 & 5 & 16  \\
  13 & 1 & 2 & 1 & \pz7              & 61 & 1 & \pz8 & 6 & 17  \\
  16 & 1 & 3 & 1 & \pz8              & 70 & 1 & \pz9 & 7 & 18  \\
  19 & 1 & 4 & 1 & \pz9              & 73 & 1 & \pz8 & 8 & 19  \\
  25 & 1 & 2 & 3 &   10              & 79 & 1 &   10 & 8 & 20  \\
\bottomrule
\end{tabular}
\caption{Maximum polynomial degree $\maxpoly$ per number of points $\npoints$.}
\vspace*{-.75em}
\label{tab:polynomial_degrees}
\end{table}

\reviewerOne{For an arbitrary function sequence, for each $\npoints$, we construct a function sequence $\mathbf{f}(x,y)=\mathbf{f}_{\ngroup}(x,y)=\{\mathbf{f}^0(x,y),\hdots,\mathbf{f}^\ngroup(x,y)\}$, using the same value for $\ngroup$ as the polynomial rules.  In doing so, we reduce the maximum polynomial degree $\maxpoly$ that can be integrated exactly from Table~\ref{tab:polynomial_degrees} in exchange for the ability to integrate the singular functions.}

\reviewerOne{For example, for $\npoints=6$, $\ngroup=4$ and $\mathbf{f}(x,y)=\{\mathbf{f}^0(x,y),\hdots,\mathbf{f}^4(x,y)\}$.  For the polynomial rules, 
$\mathbf{f}^0(x,y) = \{1\}$, 
$\mathbf{f}^1(x,y) = \{x\}$, 
$\mathbf{f}^2(x,y)=\{x^2,\allowbreak\,xy\}$,
$\mathbf{f}^3(x,y)=\{x^3,\allowbreak\,x^2y\}$, and
$\mathbf{f}^4(x,y)=\{x^4,\allowbreak\,x^3y,\allowbreak\,x^2y^2\}$.
%
If the arbitrary function sequence has two singularities, $\singfuncvec(x,y)=\{s_1(x,y),\,s_2(x,y)\}$, the functions in the function groups with the highest polynomial degrees are replaced by the singularities, such that $\mathbf{f}^3(x,y)=\{s_1(x,y)\}$ and $\mathbf{f}^4(x,y)=\{s_2(x,y)\}$.  The maximum polynomial degree is reduced to $\maxpoly=2$ and, from Section~\ref{sec:singularity_2d},~$\nfunc=6$.}

When constructing the sequence of functions, one must weigh the amount of singular functions against the maximum polynomial degree that can be integrated.  Additionally, whereas the ability to integrate polynomials includes the ability to integrate cross terms (e.g., the ability to integrate $x^3$ also indicates the ability to integrate $x^2 y$), the ability to integrate singular functions does not extend to cross terms.  Therefore, there are three approaches to address this issue:
\begin{enumerate}
\item Use a two-dimensional characterization of the singularities (Section~\ref{sec:singularity_2d}), if it is available.
\item Use a one-dimensional characterization of the singularities (Section~\ref{sec:singularity_1d}), assuming the cross terms do not warrant an additional reduction in the maximum polynomial degree that can be integrated.
\item Include cross terms for the one-dimensional characterization (Section~\ref{sec:singularity_1d}) at the expense of reducing the maximum polynomial degree that can be integrated.
\end{enumerate}
\vspace*{-.25em}
Alternatively, one can use Approach 2, which is presented in Section~\ref{sec:large}.

\reviewerOne{With the function sequence constructed, we solve the optimization problem in Equation~\eqref{eq:opt_bary} with the expectation that the objective function in Equation~\eqref{eq:objective_function} is zero.  Upon doing so, we attempt to include subsequent groups of functions until the objective function is nonzero.  Additionally, we reject points that are not within the interior of the triangle.}

\reviewerOne{%
Equation~\eqref{eq:opt_bary} is a nonlinear least squares problem, which we solve using the Levenberg--Marquardt~\cite{levenberg_1944,marquardt_1963} implementation in the \texttt{lmder} subroutine of \texttt{minpack}~\cite{minpack}.  We use the arbitrary precision library \texttt{MPFUN2015}~\cite{bailey_2019} with 64 digits of working precision and require the objective function in Equation~\eqref{eq:objective_function} to be less than $10^{-150}$.}

As mentioned in the introduction, the ability to compute quadrature rules is heavily dependent upon the initial guess for the iterative solver.  \reviewerOne{Therefore, we use points and weights close to those presented in Reference~\cite{dunavant_1985} for polynomials as initial guesses.  Additional details are provided in Appendix~\ref{app:initial}.  Additionally, the rules may not be unique; we provide a method to address this shortcoming in Appendix~\ref{app:unique}.}

From barycentric coordinates, the points can be mapped to an arbitrary triangle with vertices $(x_1,y_1)$, $(x_2,y_2)$, and $(x_3,y_3)$ by
\begin{align*}
x = \alpha x_1 + \beta x_2 + (1-\alpha-\beta) x_3, \qquad y = \alpha y_1 + \beta y_2 + (1-\alpha-\beta) y_3.
\end{align*}
The weights are multiplied by twice the area of the arbitrary triangle: $w=2Aw'$.

\section{Approach 2: Quadrilateral Subdomains}
\label{sec:large}

For polynomials, symmetric quadrature rules have been computed for high degrees.  This is facilitated by the equation reduction achieved by exploiting invariance, as mentioned in Section~\ref{sec:polynomials}.  Nonetheless, for higher degrees, the number of possible solutions increases~\cite{papanicolopulos_2015,papanicolopulos_2015_anc}.  Additionally, for higher degrees, the optimal number of points is less straightforward.

For arbitrary function sequences, which are the focus of this paper, a systematic reduction in the number of equations is unavailable, the solution to Equation~\eqref{eq:opt_bary} is more susceptible to nonzero local minima, and knowledge of the optimal number of points is unavailable.  Therefore, for large amounts of functions, we employ $\npoints'$-point one-dimensional rules that integrate the one-dimensional function sequences described in Section~\ref{sec:singularity_1d}.  
Because the rules are one-dimensional, $\nfunc=2\npoints'=\nmono+\nsing$\reviewerTwo{, where $\npoints'$ denotes the number of points in one dimension}.

To achieve symmetric quadrature rules, we do the following:
\begin{enumerate}
\item Compute one-dimensional rules $\xi_i$, $w_i'$, for $i=1,\hdots,\npoints'$ on the unit interval $\xi\in[0,\,1]$ using the continuation approach presented in Reference~\cite{ma_1996}.
\item Take the outer product of the one-dimensional rules to obtain rules for the unit square $(\xi,\eta)\in[0,\,1]\times[0,\,1]$, such that $(\xi_i,\eta_j)=(\xi_i,\xi_j)$ and $w_{ij}'=w_i' w_j'$.
\item Using the vertices, edge midpoints, and centroid, subdivide the triangle into 3 quadrilaterals: $(A,D,O,F)$, $(B,E,O,D)$, and $(C,F,O,E)$, where $O$ denotes the centroid.
\item Bilinearly transform the unit square to each quadrilateral:
\begin{align*}
x(\xi,\eta) = \sum_{k=1}^n x_k\psi_k(\xi,\eta), \qquad y(\xi,\eta) = \sum_{k=1}^n y_k\psi_k(\xi,\eta),
\end{align*}
where
\begin{align*}
\psi_1(\xi,\eta) = (1-\xi)(1-\eta), \qquad
\psi_2(\xi,\eta) = \xi(1-\eta), \qquad
\psi_3(\xi,\eta) = \xi\eta, \qquad
\psi_4(\xi,\eta) = (1-\xi)\eta, 
\end{align*}
and $\{x_1,\,x_2,\,x_3,\,x_4\}=\{x_A,\,x_D,\,x_O,\,x_F\}$, $\{x_B,\,x_E,\,x_O,\,x_D\}$, or $\{x_C,\,x_F,\,x_O,\,x_E\}$, and similarly for $y_k$.
\item Compute the weights $w_{ij}=\reviewerOne{|J(\xi_i,\eta_j)|}w'_{ij}$, where $J(\xi_i,\eta_j)$ is the determinant of the Jacobian:
\begin{align*}
J(\xi_i,\eta_j) = \left|
\begin{matrix} 
\displaystyle\frac{\partial x}{\partial \xi} &
\displaystyle\frac{\partial x}{\partial \eta} \\[1em]
\displaystyle\frac{\partial y}{\partial \xi} &
\displaystyle\frac{\partial y}{\partial \eta} 
\end{matrix}
\right|.
\end{align*}
\end{enumerate}

\begin{figure}[b]
\centering
\input{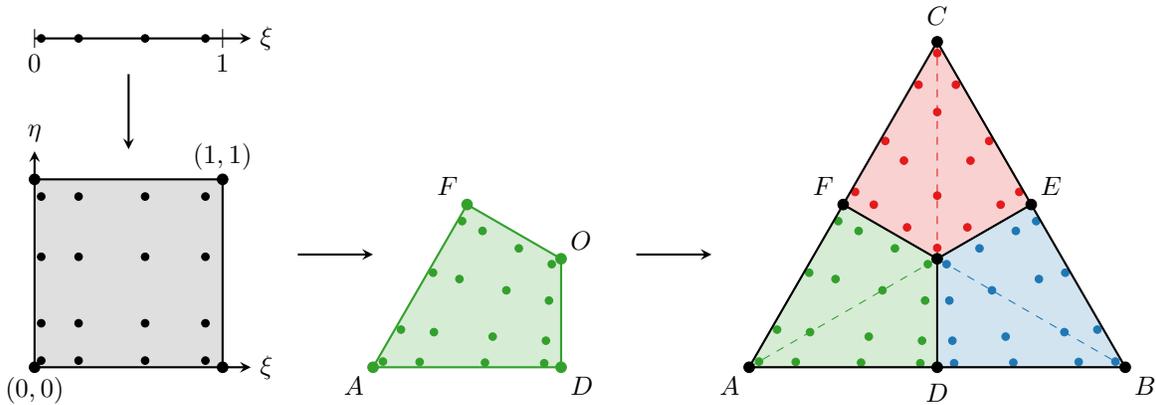}
\caption{Pictorial representation of Approach 2: A triangle is split into three quadrilateral subdomains.}
\label{fig:three_quads}
\end{figure}

The transformation is shown in Figure~\ref{fig:three_quads}.  For an $\npoints'$-point one-dimensional rule, the triangle has $\npoints=3\npoints'^2$ integration points.  Instead of performing the transformation for each triangular integration domain, the transformation can be performed on an arbitrary reference triangle, which is then linearly mapped to each triangular integration domain.  The advantage of this approach is that it is also directly applicable to quadrilaterals exhibiting boundary singularities.

\section{Numerical Example: Singular Functions Appearing in Electromagnetic Simulations}
\label{sec:example}

To demonstrate the approaches of the preceding sections, we consider the integrals
\begin{align}
I_c=\int_A \int_{A'} \frac{\cos(2\pi \|\mathbf{x}-\mathbf{x}'\|_2)}{\|\mathbf{x}-\mathbf{x}'\|_2} dA' dA, \qquad I_s=\int_A \int_{A'} \frac{\sin(2\pi \|\mathbf{x}-\mathbf{x}'\|_2)}{\|\mathbf{x}-\mathbf{x}'\|_2} dA' dA,
\label{eq:integrals}
\end{align}
which commonly appear in the integration of the free-space Green's function, $G(\mathbf{x},\mathbf{x}')\propto e^{-i\mathbf{k}\cdot(\mathbf{x}-\mathbf{x}')}/\|\mathbf{x}-\mathbf{x}'\|_2$, where $\mathbf{k}$ is real, in electromagnetic simulations.  Note that these approaches are also valid for the general case of complex $\mathbf{k}$, as reported in Reference~\cite{freno_ieee}. 
The outer integration domain $A$ is a triangle.  For the inner integration domain $A'$, we consider two domains: (1) $A'=A$ and (2) $A'$ is the co-planar reflection of $A$ along a shared edge.

For both integrals and domains, the integral over $A'$ can be integrated using a radial--angular transformation~\cite{khayat_2008,freno_ieee}.  While the integrand of $I_s$ is smooth, the integrand of the outer integral of $I_c$ exhibits logarithmic singularities along shared edges and vertices.

\subsection{\reviewerTwo{Approach 1 Function Sequences}} 
\label{sec:approach_1_function_sequence}

\reviewerTwo{We have found that a suitable choice for a one-dimensional characterization of the singularities associated with this integrand is $\singfunc_j(x)=x^{2j-1}\ln x$ for $j\in\mathbb{N}^+$~\cite{freno_ieee}, and a suitable choice for a two-dimensional characterization is }
\begin{align*}
\singfunc_j(x,y)=\left\{\begin{array}{l l}
x^j\ln\left(y-1+\sqrt{x^2+(y-1)^2}\right), &\text{for $j$ odd},\\
x^{j-1}\ln\left(y+\sqrt{x^2+y^2}\right), &\text{for $j$ even},
\end{array}\right.
\end{align*}
for $j\in\mathbb{N}^+$~\cite{freno_ieee}.

\reviewerTwo{For both characterizations, the function sequences are constructed by placing each singular function after the polynomial with the degree that matches the power of $x$ in the singular function.  The function sequences containing the one- and two-dimensional characterizations of the singularities are reported in Tables~\ref{tab:function_sequence_approach_1_1d} and~\ref{tab:function_sequence_approach_1_2d}, respectively.}  \reviewerOne{The resulting function sequence is $\mathbf{f}(x,y)=\mathbf{f}_{\ngroup}(x,y)=\{\mathbf{f}^0(x,y),\hdots,\mathbf{f}^\ngroup(x,y)\}$, and the size of the sequence is $\nfunc$.}

\begin{table}[htbp!]
\centering
\reviewerTwo{
\begin{tabular}{c c c c c}
\toprule
$\ngroup$ & $\nfunc$ & $\maxpoly$ & $\nsing$ & $\phantom{^\ngroup}\mathbf{f}^\ngroup$                                \\
\midrule       
\pz0      & \pz1     & \pz0       & \pz0     & $1$                                                 \\[.25em]                                                            
\pz1      & \pz2     & \pz1       & \pz0     & $x$                                                 \\[.25em]       
\pz2      & \pz3     & \pz1       & \pz1     & $x\ln x$                                            \\[.25em]       
\pz3      & \pz5     & \pz2       & \pz1     & $x^2,\,xy$                                          \\[.25em]       
\pz4      & \pz7     & \pz3       & \pz1     & $x^3,\,x^2y$                                        \\[.25em]       
\pz5      & \pz8     & \pz3       & \pz2     & $x^3\ln x$                                          \\[.25em]     
\pz6      &   11     & \pz4       & \pz2     & $x^4,\,x^3y,\,x^2y^2$                               \\[.25em]  
\pz7      &   14     & \pz5       & \pz2     & $x^5,\,x^4y,\,x^3y^2$                               \\[.25em]  
\pz8      &   15     & \pz5       & \pz3     & $x^5\ln x$                                          \\[.25em] 
\pz9      &   19     & \pz6       & \pz3     & $x^6,\,x^5y,\,x^4y^2,\,x^3y^3$                      \\[.25em]  
10        &   23     & \pz7       & \pz3     & $x^7,\,x^6y,\,x^5y^2,\,x^4y^3$                      \\[.25em]  
11        &   24     & \pz7       & \pz4     & $x^7\ln x$                                          \\[.25em]                                                                  
12        &   29     & \pz8       & \pz4     & $x^8,\,x^7y,\,x^6y^2,\,x^5y^3,\,x^4y^4$             \\[.25em]  
13        &   34     & \pz9       & \pz4     & $x^9,\,x^8y,\,x^7y^2,\,x^6y^3,\,x^5y^4$             \\[.25em]  
14        &   35     & \pz9       & \pz5     & $x^9\ln x$                                          \\[.25em] 
15        &   41     &   10       & \pz5     & $x^{10},\,x^9y,\,x^8y^2,\,x^7y^3,\,x^6y^4,\,x^5y^5$ \\ 
\bottomrule
\end{tabular}
}
\caption{\reviewerTwo{Approach 1: Function sequence with one-dimensional characterization of the singularities.}}
\label{tab:function_sequence_approach_1_1d}
\end{table}

\begin{table}[htbp!]
\centering
\reviewerOne{
\begin{tabular}{c c c c c}
\toprule
$\ngroup$ & $\nfunc$ & $\maxpoly$ & $\nsing$ & $\phantom{^\ngroup}\mathbf{f}^\ngroup$                                \\
\midrule       
\pz0      & \pz1     & \pz0       & \pz0     & $1$                                                 \\[.25em]                                                            
\pz1      & \pz2     & \pz1       & \pz0     & $x$                                                 \\[.25em]       
\pz2      & \pz3     & \pz1       & \pz1     & $x\ln\bigl(y-1+\sqrt{x^2+(y-1)^2}\bigr)$            \\[.25em]       
\pz3      & \pz4     & \pz1       & \pz2     & $x\ln\bigl(y+\sqrt{x^2+y^2}\bigr)$                  \\[.25em]       
\pz4      & \pz6     & \pz2       & \pz2     & $x^2,\,xy$                                          \\[.25em]       
\pz5      & \pz8     & \pz3       & \pz2     & $x^3,\,x^2y$                                        \\[.25em]       
\pz6      & \pz9     & \pz3       & \pz3     & $x^3\ln\bigl(y-1+\sqrt{x^2+(y-1)^2}\bigr)$          \\[.25em]     
\pz7      &   10     & \pz3       & \pz4     & $x^3\ln\bigl(y+\sqrt{x^2+y^2}\bigr)$                \\[.25em]      
\pz8      &   13     & \pz4       & \pz4     & $x^4,\,x^3y,\,x^2y^2$                               \\[.25em]  
\pz9      &   16     & \pz5       & \pz4     & $x^5,\,x^4y,\,x^3y^2$                               \\[.25em]  
  10      &   17     & \pz5       & \pz5     & $x^5\ln\bigl(y-1+\sqrt{x^2+(y-1)^2}\bigr)$          \\[.25em] 
  11      &   18     & \pz5       & \pz6     & $x^5\ln\bigl(y+\sqrt{x^2+y^2}\bigr)$                \\[.25em] 
  12      &   22     & \pz6       & \pz6     & $x^6,\,x^5y,\,x^4y^2,\,x^3y^3$                      \\[.25em]  
  13      &   26     & \pz7       & \pz6     & $x^7,\,x^6y,\,x^5y^2,\,x^4y^3$                      \\[.25em]  
  14      &   27     & \pz7       & \pz7     & $x^7\ln\bigl(y-1+\sqrt{x^2+(y-1)^2}\bigr)$          \\[.25em] 
  15      &   28     & \pz7       & \pz8     & $x^7\ln\bigl(y+\sqrt{x^2+y^2}\bigr)$                \\[.25em]                                                                  
  16      &   33     & \pz8       & \pz8     & $x^8,\,x^7y,\,x^6y^2,\,x^5y^3,\,x^4y^4$             \\[.25em]  
  17      &   38     & \pz9       & \pz8     & $x^9,\,x^8y,\,x^7y^2,\,x^6y^3,\,x^5y^4$             \\[.25em]  
  18      &   39     & \pz9       & \pz9     & $x^9\ln\bigl(y-1+\sqrt{x^2+(y-1)^2}\bigr)$          \\[.25em] 
  19      &   40     & \pz9       &   10     & $x^9\ln\bigl(y+\sqrt{x^2+y^2}\bigr)$                \\[.25em]
  20      &   46     &   10       &   10     & $x^{10},\,x^9y,\,x^8y^2,\,x^7y^3,\,x^6y^4,\,x^5y^5$ \\ 
\bottomrule
\end{tabular}
}
\caption{\reviewerTwo{Approach 1: Function sequence with two-dimensional characterization of the singularities.}}
\label{tab:function_sequence_approach_1_2d}
\end{table}

\reviewerOne{For these functions, we map $x\to \alpha$ and $y\to \beta$, where $\alpha$ and $\beta$ are the barycentric coordinates of a triangle.  Because the rules are geometrically symmetric, these quadrature rules are able to account for the singularities at each edge and vertex.}

Using Approach 1 from Section~\ref{sec:moderate}, for each $\npoints$, we attempt to increase $\ngroup$ \reviewerTwo{for both function sequences}.  The final values of $\ngroup$ are listed in \reviewerTwo{Tables~\ref{tab:final_ngroup_1d} and \ref{tab:final_ngroup_2d}, which provide} a comparison with the initial $\ngroup$, as well as a comparison with the maximum polynomial degree $\maxpoly$ from the polynomial rules.  

\reviewerOne{For example, as indicated in Table~\ref{tab:polynomial_degrees}, for $\npoints=16$, the polynomial rules can integrate up to $\ngroup=8$, which accounts for polynomials up to $\maxpoly=8$.  For the function sequence listed in Table~\ref{tab:function_sequence_approach_1_2d}, we initially attempt to integrate up to $\ngroup=8$, which contains polynomials up to $\maxpoly=4$ and $\nsing=4$ singularities.  These details are listed under the `Initial' block column of Table~\ref{tab:final_ngroup_2d}.  We are successfully able to integrate the function sequence associated with $\ngroup=8$; we are also able to successfully integrate the function sequences associated with $\ngroup=9$ and $\ngroup=10$.  The function sequence associated with $\ngroup=10$ includes polynomials up to $\maxpoly=5$ and $\nsing=5$ singularities.  These details are listed under the `Final' block column of Table~\ref{tab:final_ngroup_2d}.}

\begin{figure}[htbp!]
\centering
\input{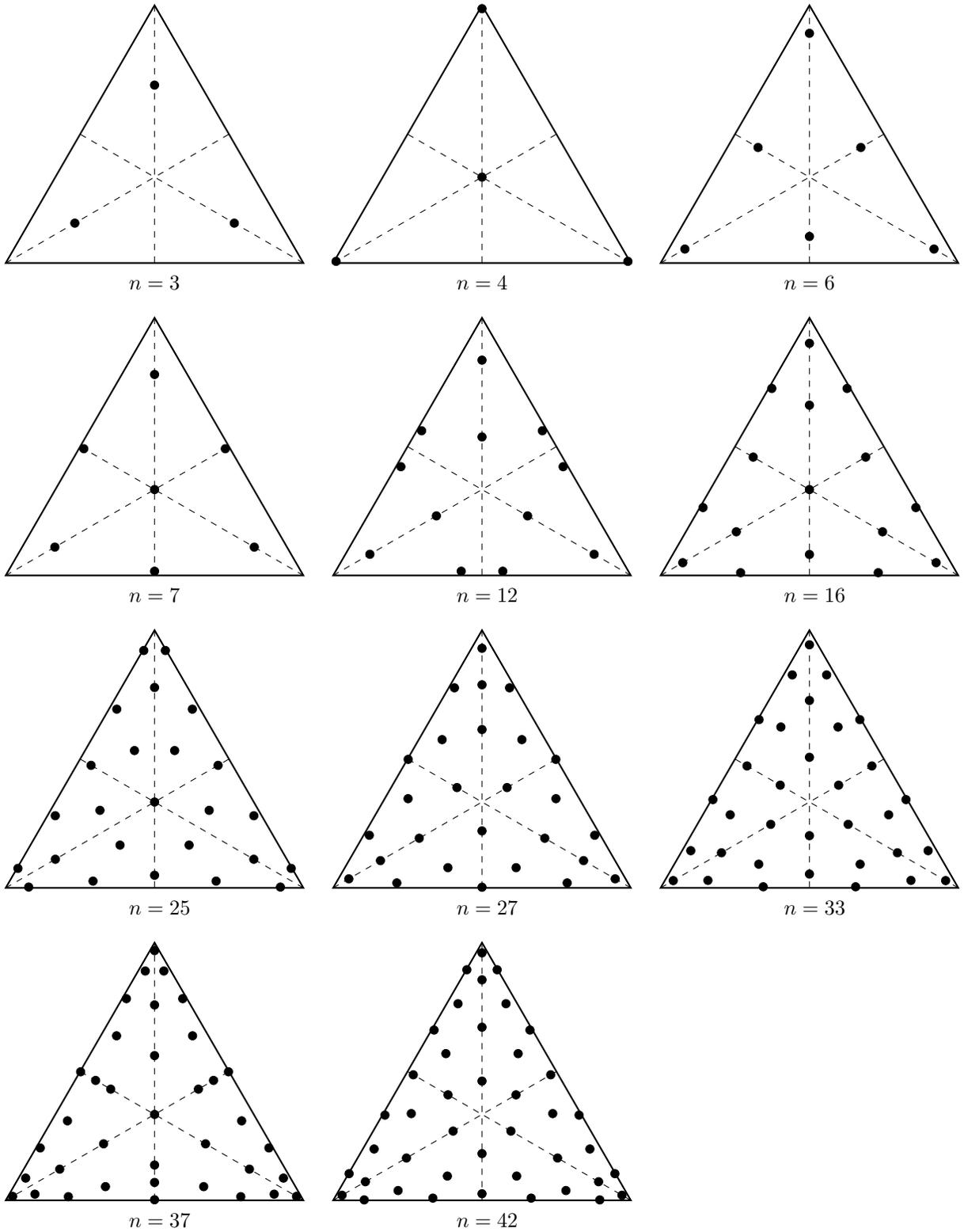}
\vspace*{-1em}
\caption{\reviewerTwo{Approach 1: Points for $\npoints$ between 3 and 42 using one-dimensional characterization of the singularities.}}
\label{fig:points_1_1d}
\end{figure}

When, for a given $\npoints$, the final $\ngroup$ is not greater than that from a lower $\npoints$, the higher $\npoints$ is eliminated.  For example,
\reviewerTwo{in Tables~\ref{tab:final_ngroup_1d} and \ref{tab:final_ngroup_2d},} for $\npoints=13$, $\ngroup$ does not increase; however, for $\npoints=12$, $\ngroup$ is increased by one, matching the final $\ngroup$ for $\npoints=13$.  Therefore, we eliminate $\npoints=13$.  
Diagrams of the points computed using Approach 1 are shown in \reviewerTwo{Figure~\ref{fig:points_1_1d} using the one-dimensional characterization of the singularities and in Figure~\ref{fig:points_1_2d} using the two-dimensional characterization of the singularities}.

\begin{table}[!h]
\centering
\reviewerTwo{
\begin{tabular}{c c c c c c c c}
\toprule
& \multicolumn{3}{c}{Initial} & \multicolumn{3}{c}{Final} & \multicolumn{1}{c}{Polynomial Rules} \\
\cmidrule(lr){2-4} \cmidrule(lr){5-7} \cmidrule(lr){8-8}
$\npoints$ & $\ngroup$ & $\maxpoly$ & $\nsing$ & $\ngroup$ & $\maxpoly$ & $\nsing$ & $\ngroup=\maxpoly$  \\
\midrule
\pz1       & \pz1      & \pz1       & 0        & \pz1      & \pz1       & \pz0     & \pz1 \\
\pz3       & \pz2      & \pz1       & 1        & \pz2      & \pz1       & \pz1     & \pz2 \\
\pz4       & \pz3      & \pz2       & 1        & \pz3      & \pz2       & \pz1     & \pz3 \\
\pz6       & \pz4      & \pz3       & 1        & \pz4      & \pz3       & \pz1     & \pz4 \\
\pz7       & \pz5      & \pz3       & 2        & \pz5      & \pz3       & \pz2     & \pz5 \\
  12       & \pz6      & \pz4       & 2        & \pz7      & \pz5       & \pz2     & \pz6 \\
  13       & \pz7      & \pz5       & 2        & $-$       & $-$        & $-$      & \pz7 \\
  16       & \pz8      & \pz5       & 3        & \pz8      & \pz5       & \pz3     & \pz8 \\
  19       & \pz9      & \pz6       & 3        & $-$       & $-$        & $-$      & \pz9 \\
  25       &   10      & \pz7       & 3        &   12      & \pz8       & \pz4     &   10 \\
  27       &   11      & \pz7       & 4        &   13      & \pz9       & \pz4     &   11 \\
  33       &   12      & \pz8       & 4        &   15      &   10       & \pz5     &   12 \\
  37       &   13      & \pz9       & 4        &   16      &   11       & \pz5     &   13 \\
  42       &   14      & \pz9       & 5        &   17      &   11       & \pz6     &   14 \\
\bottomrule
\end{tabular}
}
\caption{\reviewerTwo{Approach 1: Function groups integrated with one-dimensional characterization of the singularities.}}
\label{tab:final_ngroup_1d}
\end{table}

\begin{table}[!h]
\centering
\begin{tabular}{c c c c c c c c}
\toprule
& \multicolumn{3}{c}{Initial} & \multicolumn{3}{c}{Final} & \multicolumn{1}{c}{Polynomial Rules} \\
\cmidrule(lr){2-4} \cmidrule(lr){5-7} \cmidrule(lr){8-8}
$\npoints$ & $\ngroup$ & $\maxpoly$ & $\nsing$ & $\ngroup$ & $\maxpoly$ & $\nsing$ & $\ngroup=\maxpoly$  \\
\midrule
\pz1       & \pz1      &  1         & 0        & \pz1      & \pz1          & \pz0        & \pz1 \\
\pz3       & \pz2      &  1         & 1        & \pz2      & \pz1          & \pz1        & \pz2 \\
\pz4       & \pz3      &  1         & 2        & \pz3      & \pz1          & \pz2        & \pz3 \\
\pz6       & \pz4      &  2         & 2        & \pz4      & \pz2          & \pz2        & \pz4 \\
\pz7       & \pz5      &  3         & 2        & \pz5      & \pz3          & \pz2        & \pz5 \\
  12       & \pz6      &  3         & 3        & \pz7      & \pz3          & \pz4        & \pz6 \\
  13       & \pz7      &  3         & 4        & $-$       & $-$           & $-$         & \pz7 \\
  16       & \pz8      &  4         & 4        &   10      & \pz5          & \pz5        & \pz8 \\
  19       & \pz9      &  5         & 4        &   11      & \pz5          & \pz6        & \pz9 \\
  25       &   10      &  5         & 5        &   13      & \pz7          & \pz6        &   10 \\
  27       &   11      &  5         & 6        &   15      & \pz7          & \pz8        &   11 \\
  33       &   12      &  6         & 6        &   16      & \pz8          & \pz8        &   12 \\
  37       &   13      &  7         & 6        & $-$       & $-$           & $-$         &   13 \\
  42       &   14      &  7         & 7        &   20      & 10            & 10          &   14 \\
  48       &   15      &  7         & 8        & $-$       & $-$           & $-$         &   15 \\
  52       &   16      &  8         & 8        &   23      & 11            & 12          &   16 \\
\bottomrule
\end{tabular}
\caption{\reviewerTwo{Approach 1: Function groups integrated with two-dimensional characterization of the singularities.}}
\label{tab:final_ngroup_2d}
\end{table}

Though not shown in this paper, some missing final $\ngroup$ values, including those for eliminated points, can be obtained through the following approaches (Reference~\cite{papanicolopulos_2015_anc} lists many of these options):
\begin{enumerate}
\item For a given $\npoints$ and orbit amount $(n_0,n_1,n_2)$, use \reviewerOne{initial guesses near another set of polynomial rules}.  As stated previously, the polynomial rules are often not unique and the solution to Equation~\eqref{eq:opt_bary} is dependent upon initial guess.
\item For a given $\npoints$, use a different polynomial-rule orbit amount.  For example, for $n=37$, Reference~\cite{dunavant_1985} provides a $(1,6,3)$ rule, whereas Reference~\cite{papanicolopulos_2015} provides a $(1,4,4)$ rule.
\item Use a slightly higher $\npoints$ that is less efficient as a polynomial rule.  For example, a 13-point $(1,2,1)$ rule and a 15-point $(0,1,2)$ rule can both integrate polynomials of degree $\maxpoly=7$~\cite{papanicolopulos_2015}; however, \reviewerTwo{for the two-dimensional characterization of the singularities}, using the 13-point rule from Reference~\cite{dunavant_1985} yields a final $\ngroup=7$, whereas one of the 15-point rules from Reference~\cite{papanicolopulos_2015_anc} yields a final $\ngroup=8$.
\end{enumerate}

\begin{figure}[htbp!]
\centering
\input{points_1_2d.tex}
\vspace*{-1em}
\caption{\reviewerTwo{Approach 1: Points for $\npoints$ between 3 and 52 using two-dimensional characterization of the singularities.}}
\label{fig:points_1_2d}
\end{figure}

\newpage

\subsection{\reviewerTwo{Approach 2 Function Sequences}} 

\reviewerTwo{Using the same one-dimensional characterization of the singularities as for Approach 1 in Section~\ref{sec:approach_1_function_sequence}, the one-dimensional function sequence for Approach 2 is reported in Table~\ref{tab:function_sequence_approach_2}.}

\begin{table}[htbp!]
\centering
\begin{tabular}{c c c}
\toprule
$\npoints$ & $\npoints'$ & Functions \\
\midrule      
\pz\pz3 & 1 & $1,\,x$          \\             
  \pz12 & 2 & $x\ln x,\,x^2$   \\     
  \pz27 & 3 & $x^3,\,x^3\ln x$ \\ 
  \pz48 & 4 & $x^4,\,x^5$      \\
  \pz75 & 5 & $x^5\ln x,\,x^6$ \\     
    108 & 6 & $x^7,\,x^7\ln x$ \\ 
\bottomrule
\end{tabular}
\caption{\reviewerTwo{Approach 2: Function sequence.}}
\label{tab:function_sequence_approach_2}
\end{table}

\begin{figure}[htbp!]
\centering
\input{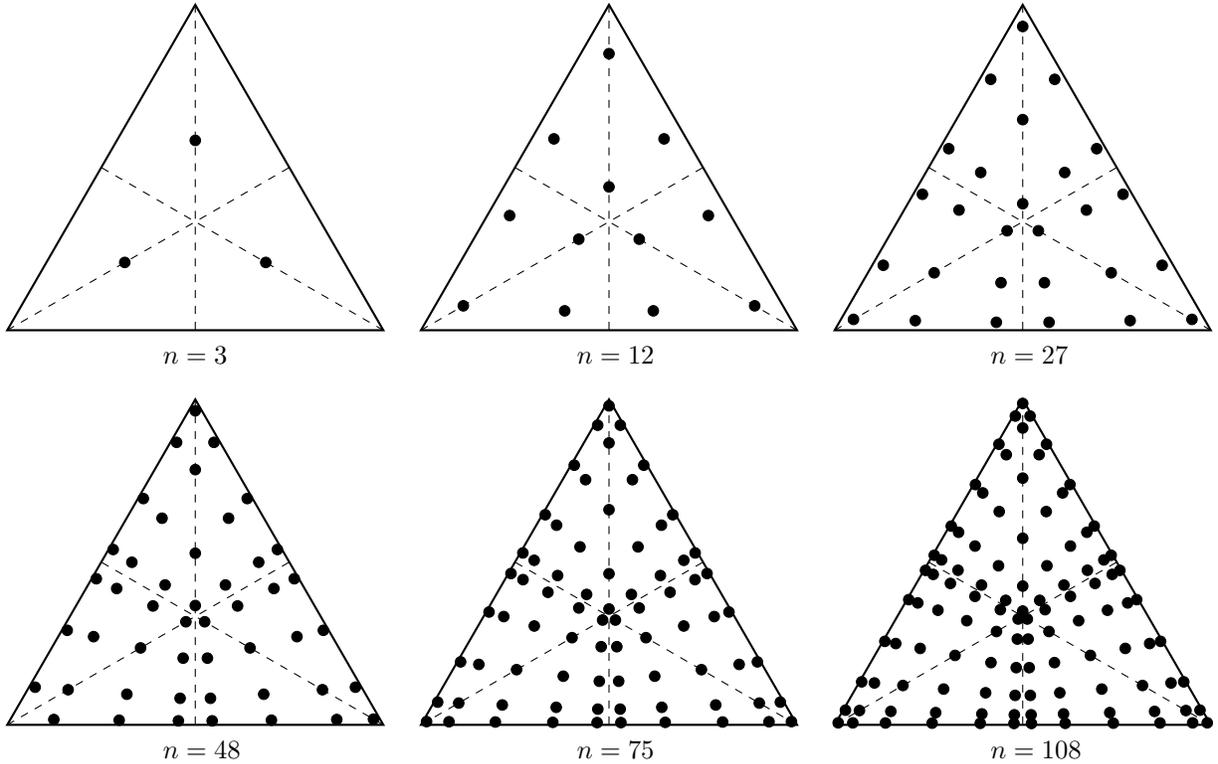}
\vspace*{-1em}
\caption{\reviewerTwo{Approach 2: Points $\npoints$ between 3 and 108.}}
\label{fig:points_2}
\end{figure}

Using the functions in Table~\ref{tab:function_sequence_approach_2}, we construct the one-dimensional rules used by Approach 2.  Diagrams of the points computed using Approach 2 are shown in Figure~\ref{fig:points_2}.

\subsection{Results} 
To assess the performance of the rules arising from these approaches, we compute reference solutions for $I_c$ and $I_s$, as described in Reference~\cite{freno_ieee}, for $A$ being defined by the vertices $(0,0)$, $(1/20,1/20)$, and $(-1/20,1/20)$.  For Domain 1, $A'=A$ and, for Domain 2, $A'$ is defined by the vertices $(0,1/10)$, $(-1/20,1/20)$, and $(1/20,1/20)$.

\begin{figure}
\centering
\begin{subfigure}[b]{.49\textwidth}
\begin{tikzpicture}
\node at (0,0) {\includegraphics[scale=.6,clip=true,trim=2.35in 0in 2.87in 0in]{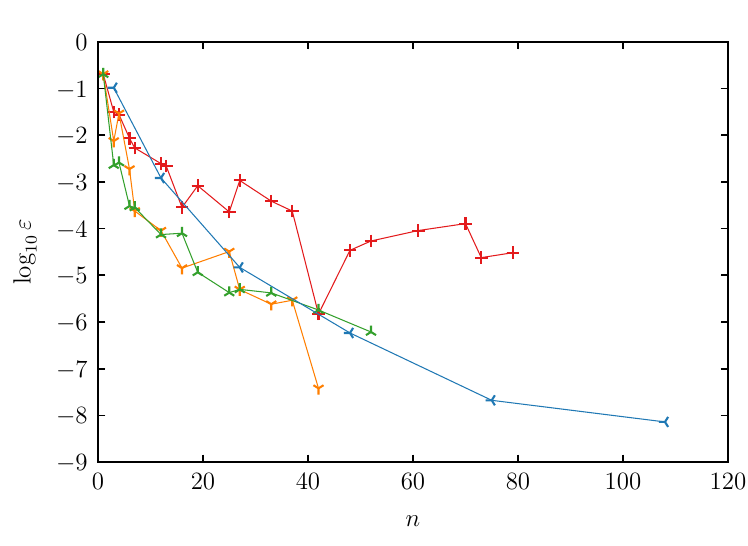}};
\def\ts{15}; 
\def\tx{2.6}; 
\def\ty{1.4}; 
\coordinate (A) at ($\ts*( 0.00,0.00)+(\tx,\ty)$);
\coordinate (B) at ($\ts*( 0.05,0.05)+(\tx,\ty)$);
\coordinate (C) at ($\ts*(-0.05,0.05)+(\tx,\ty)$);
\coordinate (E) at ($1/3*(A)+1/3*(B)+1/3*(C)$);
\draw[line cap=round,line join=round] (A) -- (B) -- (C) -- cycle;
\node at (E) {$A,A'$};
\end{tikzpicture}
\caption{$I_c$}
\label{fig:cos_180}
\end{subfigure}
~ 
\begin{subfigure}[b]{.49\textwidth}
\begin{tikzpicture}
\node at (0,0) {\includegraphics[scale=.6,clip=true,trim=2.35in 0in 2.87in 0in]{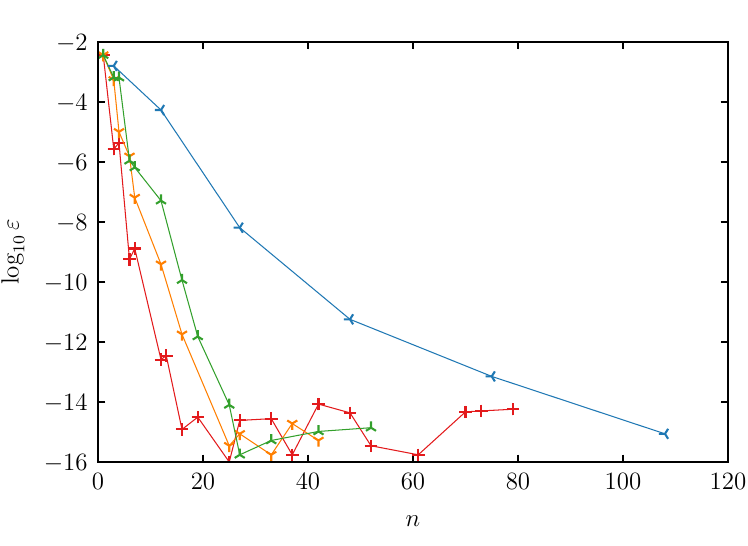}};

\node at (1.5,1.3) {
\includegraphics[scale=.6,clip=true,trim=7.3in 1.4in 0.4in 1.1in]{t1_180_sin.pdf}}
;
\end{tikzpicture}
\caption{$I_s$}
\label{fig:sin_180}
\end{subfigure}
\caption{Domain 1: Relative errors for $I_c$ and $I_s$ when $A'=A$.}
\label{fig:180}
\end{figure}

\begin{figure}
\centering
\begin{subfigure}[b]{.49\textwidth}
\begin{tikzpicture}
\node at (0,0) {\includegraphics[scale=.6,clip=true,trim=2.35in 0in 2.87in 0in]{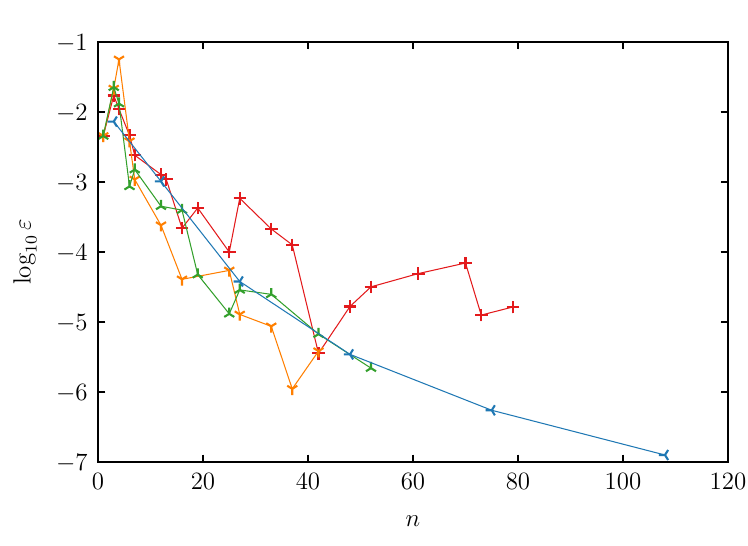}};
\def\ts{15}; 
\def\tx{2.6}; 
\def\ty{.65}; 
\coordinate (A) at ($\ts*( 0.00,0.00)+(\tx,\ty)$);
\coordinate (B) at ($\ts*( 0.05,0.05)+(\tx,\ty)$);
\coordinate (C) at ($\ts*(-0.05,0.05)+(\tx,\ty)$);
\coordinate (D) at ($\ts*( 0.00,0.10)+(\tx,\ty)$);
\coordinate (E) at ($1/3*(A)+1/3*(B)+1/3*(C)$);
\coordinate (F) at ($1/3*(D)+1/3*(B)+1/3*(C)$);
\draw[line cap=round,line join=round] (A) -- (B) -- (C) -- cycle;
\draw[line cap=round,line join=round] (D) -- (B) -- (C) -- cycle;

\node at (E) {$A$};
\node at (F) {$A'$};
\end{tikzpicture}
\caption{$I_c$}
\label{fig:cos_0}
\end{subfigure}
~ 
\begin{subfigure}[b]{.49\textwidth}
\begin{tikzpicture}
\node at (0,0) {\includegraphics[scale=.6,clip=true,trim=2.35in 0in 2.87in 0in]{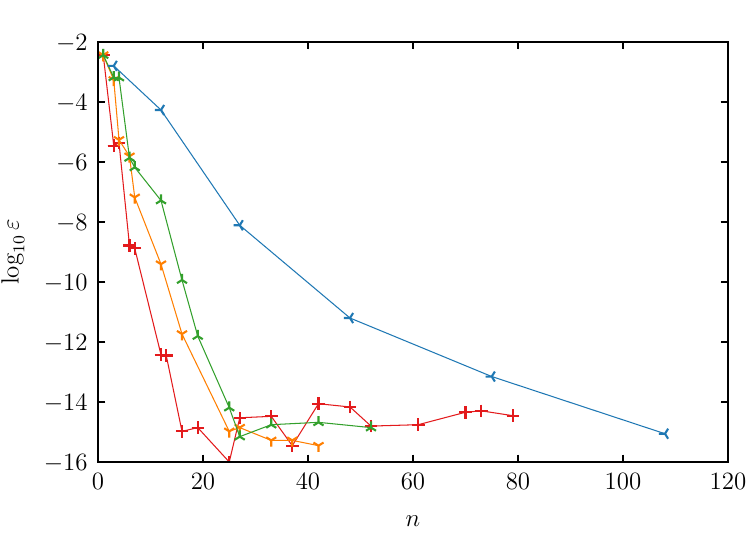}};

\node at (1.5,1.3) {
\includegraphics[scale=.6,clip=true,trim=7.3in 1.4in 0.4in 1.1in]{t1_0_sin.pdf}}
;
\end{tikzpicture}
\caption{$I_s$}
\label{fig:sin_0}
\end{subfigure}
\caption{Domain 2: Relative errors for $I_c$ and $I_s$ when $A'\ne A$.}
\label{fig:0}
\end{figure}

Figures~\ref{fig:180}--\ref{fig:0} show the relative errors $\varepsilon=\left|(\tilde{I}_k-I_k)/I_k\right|$, between the reference solution $I_k$ and the quadrature solution $\tilde{I}_k$ for $k=\{c,\,s\}$ in Equation~\ref{eq:integrals}.  These figures compare the quadrature rules presented in this paper, as well as the polynomial quadrature rules for both domains and integrals.
As shown in Figures~\ref{fig:cos_180} and~\ref{fig:cos_0} for $I_c$, both approaches generally outperform the polynomial quadrature rules, and Approach 1 often outperforms the polynomial quadrature rules by orders of magnitude.  \reviewerTwo{In Figure~\ref{fig:cos_180}, for example, Approach 1 outperforms the polynomial rules by two orders of magnitude for $\npoints=27$}.  Approach 1 is then appropriate when a moderate accuracy (e.g., 6 or 7 digits) is sufficient.  \reviewerTwo{Though the two-dimensional characterization of the singularities more accurately describes them, the appearance of the singularities in the function sequence in pairs results in a greater reduction in the maximum polynomial degree that can be integrated than the reduction due to the one-dimensional characterization.}

Because the integrand for $I_s$ is not singular, the polynomial rules perform the best; however, Approach~1 has similar efficiency, as shown in Figures~\ref{fig:sin_180} and~\ref{fig:sin_0}.

Though not as efficient as Approach 1 for $I_c$ or the polynomial rules for $I_s$, the relative error arising from Approach 2 decreases monotonically with respect to $\npoints$, which is an important feature to guarantee improved accuracy when increasing $\npoints$.  Additionally, for large values of $\npoints$, the points arising from Approach 2 take less time to compute than those from Approach 1 since they arise from one-dimensional rules.


\clearpage
\section{Conclusions and Outlook for Future Work}
\label{sec:conclusions}

In this paper, we presented two approaches for computing symmetric triangle quadrature rules for arbitrary functions.  On an example problem, we demonstrated that, for several point amounts, our approaches can achieve relative errors two orders of magnitude less than those achieved from polynomial rules. 

The more novel Approach 1 of Section~\ref{sec:moderate} yields efficient points, but, for larger amounts of points, their computation is costly and knowledge of their optimality and uniqueness is limited.  Approaches to improve these shortcomings would be beneficial, though most likely highly dependent upon the function sequence.
\section*{Acknowledgments} 
\label{sec:acknowledgments}
The authors thank Prof.\ Donald Wilton from the University of Houston for his insightful discussions.
This paper describes objective technical results and analysis. Any subjective views or opinions that might be expressed in the paper do not necessarily represent the views of the U.S. Department of Energy or the United States Government.
Sandia National Laboratories is a multimission laboratory managed and operated by National Technology and Engineering Solutions of Sandia, LLC, a wholly owned subsidiary of Honeywell International, Inc., for the U.S. Department of Energy's National Nuclear Security Administration under contract DE-NA-0003525.

\appendix
\renewcommand*{\thesection}{Appendix \Alph{section}}
\renewcommand*{\thesubsection}{\Alph{section}.\arabic{subsection}}
\section{Computational Techniques for Approach 1}

\subsection{Initial Guess}
\label{app:initial}
\reviewerOne{%
The ability to compute quadrature rules is heavily dependent upon the initial guess for the iterative solver.  Therefore, we use points and weights close to those presented in Reference~\cite{dunavant_1985} for polynomials as initial guesses.  Points and weights near those of other polynomial rules can be used, such as those listed in Reference~\cite{papanicolopulos_2015_anc}.  

Because they contribute linearly to the quadrature computation and are therefore less volatile, we use the polynomial weights directly as initial guesses.  For type-1 and type-2 orbit points, the initial guesses are sampled randomly, within one tenth of a median of the polynomial rule locations.  For the results in Section~\ref{sec:example}, we use 40 Latin hypercube samples.
}

\subsection{Uniqueness}
\label{app:unique}
\reviewerOne{%
As with the polynomial rules, for a given orbit count, the points and weights may not be unique.  To address this shortcoming, for a given $\ngroup$, we modify Equation~\eqref{eq:objective_function} to include the functions introduced by the next group, which we weight by $\theta$:
\begin{align}
F(\boldsymbol{\alpha},\boldsymbol{\beta},\mathbf{w};\theta) 
=
\sum_{f\in\mathbf{f}_{\ngroup}} \left(\frac{\tilde{I}_{f}-I_{f}}{I_{f}}\right)^2
+
\theta\!\!\!\!\sum_{f\in\mathbf{f}^{\ngroup+1}} \!\!\!\left(\frac{\tilde{I}_{f}-I_{f}}{I_{f}}\right)^2,
\label{eq:objective_function_modified}
\end{align}
where $\mathbf{f}_{\ngroup}=\{\mathbf{f}^0,\hdots,\mathbf{f}^\ngroup\}$ and $\mathbf{f}^{\ngroup+1}=\mathbf{f}_{\ngroup+1}\setminus\mathbf{f}_{\ngroup}$.

Initially, we take $\theta=1$.  Upon solving Equation~\eqref{eq:opt_bary}, Equation~\eqref{eq:objective_function_modified} is nonzero.  
We reduce $\theta$ by a factor of 10 in Equation~\eqref{eq:objective_function_modified} and compute the new solution to Equation~\eqref{eq:opt_bary} using the previous solution as an initial guess.
We repeat this process until Equation~\eqref{eq:objective_function} is zero.  When multiple solutions result in Equation~\eqref{eq:objective_function} being zero, we select the solution that yields the lowest norm of the error of the additional functions: $\displaystyle\!\!\sum_{f\in\mathbf{f}^{\ngroup+1}} \!\!\!\!\left(\frac{\tilde{I}_{f}-I_{f}}{I_{f}}\right)^2$.
}


\addcontentsline{toc}{section}{\refname}
\bibliographystyle{elsarticle-num}
\bibliography{quadrature}

\end{document}